\documentclass[preprint,3p,times]{elsarticle}

\usepackage{amssymb}
\usepackage{amsmath}
\usepackage[all,cmtip]{xy}
\usepackage{amsthm}
\usepackage{color}
\usepackage{enumitem}
\usepackage{tikz}
\usepackage{appendix}
\usepackage[colorlinks=true]{hyperref}
\usepackage[hyphenbreaks]{breakurl}
\input diagxy

\theoremstyle{plain}
  \newtheorem{thm}{Theorem}[section]
  \newtheorem{lem}[thm]{Lemma}
  \newtheorem{prop}[thm]{Proposition}
  \newtheorem{cor}[thm]{Corollary}
\theoremstyle{definition}
  \newtheorem{defn}[thm]{Definition}
  \newtheorem{exmp}[thm]{Example}
  \newtheorem{rem}[thm]{Remark}
  
  \newtheorem{ques}[thm]{Question}

\DeclareMathAlphabet{\mathcal}{OMS}{cmsy}{m}{n}

\DeclareMathOperator{\id}{id}
\DeclareMathOperator{\ob}{ob}

\DeclareMathOperator{\Int}{Int}

\makeatletter
\def\ps@pprintTitle{%
 \let\@oddhead\@empty
 \let\@evenhead\@empty
 \def\@oddfoot{\centerline{\thepage}}%
 \let\@evenfoot\@oddfoot}
\makeatother

\def\oto{{\bfig\morphism<180,0>[\mkern-4mu`\mkern-4mu;]\place(86,0)[\circ]\efig}}
\def\rto{{\bfig\morphism<180,0>[\mkern-4mu`\mkern-4mu;]\place(82,0)[\mapstochar]\efig}}

\newcommand{\ra}{\rightarrow}
\newcommand{\la}{\leftarrow}

\newcommand{\lda}{\swarrow}
\newcommand{\rda}{\searrow}
\newcommand{\bs}{\backslash}

\newcommand{\bv}{\bigvee}
\newcommand{\bw}{\bigwedge}

\newcommand{\dv}{\dashv}

\newcommand{\nat}{\natural}

\renewcommand{\phi}{\varphi}
\newcommand{\al}{\alpha}
\newcommand{\be}{\beta}

\newcommand{\ga}{\gamma}

\newcommand{\lam}{\lambda}

\newcommand{\Om}{\Omega}

\newcommand{\CA}{\mathcal{A}}

\newcommand{\CC}{\mathcal{C}}
\newcommand{\CD}{\mathcal{D}}

\newcommand{\CI}{\mathcal{I}}

\newcommand{\CO}{\mathcal{O}}

\newcommand{\CQ}{\mathcal{Q}}

\newcommand{\sG}{\mathsf{G}}

\newcommand{\sP}{\mathsf{P}}
\newcommand{\sQ}{\mathsf{Q}}

\newcommand{\sfc}{\mathsf{c}}

\newcommand{\sy}{\mathsf{y}}
\newcommand{\sfs}{\mathsf{s}}
\newcommand{\Fix}{\mathsf{Fix}}

\newcommand{\bbP}{\mathbb{P}}

\newcommand{\BD}{{\bf D}}

\newcommand{\FP}{\mathfrak{P}}

\newcommand{\FU}{\mathfrak{U}}

\newcommand{\Cat}{{\bf Cat}}

\newcommand{\Dist}{{\bf Dist}}

\newcommand{\Ord}{{\bf Ord}}

\newcommand{\Rel}{{\bf Rel}}
\newcommand{\Set}{{\bf Set}}

\newcommand{\Sup}{{\bf Sup}}

\newcommand{\QCat}{\CQ\text{-}\Cat}
\newcommand{\QSymCat}{\CQ\text{-}{\bf SymCat}}

\newcommand{\QSet}{\sQ\text{-}\Set}

\newcommand{\QDist}{\CQ\text{-}\Dist}

\newcommand{\QRel}{\CQ\text{-}\Rel}
\newcommand{\QSup}{\CQ\text{-}\Sup}

\newcommand{\sPd}{\sP^{\dag}}
\newcommand{\syd}{\sy^{\dag}}

\newcommand{\op}{{\rm op}}

\newcommand{\PX}{\sP X}
\newcommand{\PY}{\sP Y}

\newcommand{\PdX}{\sPd X}

\newcommand{\HsQ}{\BD_*(\sQ)}

\newcommand{\Ps}{\sP_{\sfs}}
\newcommand{\bbPs}{\bbP_{\sfs}}
\newcommand{\Us}{\FU_{\sfs}}

\newcommand{\OX}{\CO(X)}

\newcommand{\PCX}{\mathsf{PC}(X)}
\newcommand{\ldd}{\mathrel{/}}
\newcommand{\rdd}{\mathrel{\bs}}
\newcommand{\with}{\mathrel{\&}}

\renewcommand{\leq}{\leqslant}
\renewcommand{\geq}{\geqslant}

\numberwithin{equation}{section}

\begin{document}

\begin{frontmatter}



\title{The powerset monad on quantale-valued sets}


\author[S]{Lili Shen}
\ead{shenlili@scu.edu.cn}

\author[Z]{Xiaojuan Zhao\corref{cor}}
\ead{xjzhao@swjtu.edu.cn}

\cortext[cor]{Corresponding author.}
\address[S]{School of Mathematics, Sichuan University, Chengdu 610064, China}
\address[Z]{Department of Mathematics, Southwest Jiaotong University, Chengdu 611756, China}

\begin{abstract}
For a small involutive quantaloid $\mathcal{Q}$, it is shown that the category of separated complete $\mathcal{Q}$-categories and left adjoint $\mathcal{Q}$-functors is strictly monadic over the category of symmetric $\mathcal{Q}$-categories. In particular, the (covariant) powerset monad on the category of quantale-valued sets is precisely formulated.
\end{abstract}


\end{frontmatter}




\section{Introduction}

The (covariant) powerset monad 
\begin{equation} \label{P-monad}
\bbP=(\sP,\{-\},\cup)
\end{equation}
on the category $\Set$ is well known (see, e.g., \cite[Example II.3.1.1]{Hofmann2014} or \cite[Example 5.1.5]{Riehl2016}). Explicitly:
\begin{itemize}
\item the functor $\sP\colon\Set\to\Set$ sends each (crisp) set $X$ to its \emph{powerset} $\PX$, and each map $f\colon X\to Y$ to 
\begin{equation} \label{fra}
f^{\ra}\colon\PX\to\PY,\quad A\mapsto\{fx\mid x\in A\};
\end{equation}
\item the unit is given by
\begin{equation} \label{y}
\{-\}\colon X\to\PX,\quad x\mapsto\{x\};
\end{equation}
\item the multiplication is given by
\begin{equation} \label{bigcup}
\cup\colon\sP\PX\to\PX,\quad\CA\mapsto\bigcup\CA.
\end{equation}
\end{itemize}
The Eilenberg-Moore category of this monad is exactly the category $\Sup$ of complete lattices and join-preserving maps. In other words, $\Sup$ is \emph{strictly monadic} over $\Set$. More precisely:
\begin{itemize}
\item since \eqref{fra} is always a join-preserving map between complete lattices, there is a functor $\FP\colon\Set\to\Sup$ obtained by replacing the codomain of $\sP$ with $\Sup$;
\item $\FP$ is left adjoint to the forgetful functor $\FU\colon\Sup\to\Set$, and the induced monad on $\Set$ is \eqref{P-monad};
\item the right adjoint functor $\FU\colon\Sup\to\Set$ is strictly monadic.
\end{itemize}

Now, let $\CQ$ be a small \emph{involutive quantaloid} \cite{Rosenthal1996}. From the viewpoint of category theory, it is natural to consider the $\CQ$-enriched version of the monad \eqref{P-monad}. The following results are already known:
\begin{itemize}
\item The category $\QSup$ of separated complete $\CQ$-categories and left adjoint $\CQ$-functors is strictly monadic over the category $\QCat$ of $\CQ$-categories and $\CQ$-functors \cite{Stubbe2017}.
\item $\QSup$ is strictly monadic over the slice category $\Set/\CQ_0$, where $\CQ_0$ is the set of objects of $\CQ$ \cite{Pu2015}.
\end{itemize}
When $\CQ=\mathbf{2}$, the two-element Boolean algebra, these results reduce to the strict monadicity of $\Sup$ over $\Ord$ and $\Set$, respectively, where $\Ord$ refers to the category of preordered sets and order-preserving maps. In other words, both the forgetful functors $\FU_{\sfc}\colon\QSup\to\QCat$ and $\FU=(\QSup\to^{\FU_{\sfc}}\QCat\to^{\FU_0}\Set/\CQ_0)$ are strictly monadic.
\begin{equation} \label{QSup-QCat-SetQ0}
\bfig
\btriangle|lrb|/{@{->}@/^-1.5em/}`{@{->}@/^1.5em/}`{@{->}@/^-1.5em/}/<1000,600>[\QSup`\QCat`\ \Set/\CQ_0;\FU_{\sfc}`\FU`\FU_0]
\btriangle|rba|/<-`<-`<-^)/<1000,600>[\QSup`\QCat`\ \Set/\CQ_0;\FP_{\sfc}`\FP`]
\place(-70,350)[\mbox{\rotatebox{-90}{$\bot$}}]
\place(480,-80)[\bot]
\place(540,390)[\mbox{\rotatebox{150}{$\bot$}}]
\btriangle(2000,0)|lrb|/{@{->}@/^-1.5em/}`{@{->}@/^1.5em/}`{@{->}@/^-1.5em/}/<1000,600>[\Sup`\Ord`\mbox{\ $\Set$\quad} ;\FU_{\sfc}`\FU`\FU_0]
\btriangle(2000,0)|rba|/<-`<-`<-^)/<1000,600>[\Sup`\Ord`\mbox{\ $\Set$\quad};\FP_{\sfc}`\FP`]
\place(1930,350)[\mbox{\rotatebox{-90}{$\bot$}}]
\place(2480,-80)[\bot]
\place(2540,390)[\mbox{\rotatebox{150}{$\bot$}}]
\efig
\end{equation}

Therefore, the classical notion of \emph{powerset} may be extended to the $\CQ$-enriched version as follows:
\begin{itemize}
\item The ``power'' of a $\CQ$-category $X$ is given by its image under $\FU_{\sfc}\FP_{\sfc}$ (where $\FP_{\sfc}\dv\FU_{\sfc}$), which is precisely the $\CQ$-category $\PX$ of presheaves on $X$.
\item The ``power'' of a $\CQ_0$-typed set $X$ (i.e., a set $X$ equipped with a map $|\text{-}|\colon X\to\CQ_0$) is given by its image under $\FU\FP$ (where $\FP\dv\FU$), which is precisely the underlying $\CQ_0$-typed set of the presheaf $\CQ$-category of the \emph{discrete} $\CQ$-category $X$.
\end{itemize}

Let us look again at the special case of $\CQ={\bf 2}$. Since 
\[(\sP\colon\Set\to\Set)=(\Set\to^{\FP}\Sup\to^{\FU}\Set),\] 
there are two steps to obtain the powerset of a set $X$: 
\begin{itemize}
\item first, generate the complete lattice $\PX$ of all subsets of $X$ (ordered by inclusion ``$\subseteq$'') under the functor $\FP$; 
\item second, forget the order ``$\subseteq$'' on $\PX$ under the functor $\FU$. 
\end{itemize}
As the motivation of this paper, we point out that there is another interpretation of the second step: the discrete set $\PX$ may also be regarded as the \emph{symmetrization} of the partially ordered set $(\PX,\subseteq)$. 

However, for a general (small involutive) quantaloid $\CQ$, the symmetrization of a presheaf $\CQ$-category $\PX$ is far more complicated than the underlying $\CQ_0$-typed set of $\PX$. It is now natural to ask what happens if the node $\Set/\CQ_0$ in the first triangle of \eqref{QSup-QCat-SetQ0} is replaced by $\QSymCat$, the full subcategory of $\QCat$ consisting of symmetric $\CQ$-categories. More specifically, with $(-)_{\sfs}\colon\QCat\to\QSymCat$ denoting the symmetrization functor:

\begin{ques} \label{ques}
Is the composite functor
\begin{equation} \label{Us}
\Us=(\QSup\to^{\FU_{\sfc}}\QCat\to^{(-)_{\sfs}}\QSymCat)
\end{equation}
monadic?
\end{ques}

\begin{equation} \label{QSup-QCat-QSymCat-SetQ0}
\bfig
\btriangle|lrb|/{@{->}@/^-1.5em/}`{@{->}@/^1.5em/}`{@{->}@/^-1.5em/}/<1000,600>[\QSup`\QCat`\ \QSymCat;\FU_{\sfc}`\Us`(-)_{\sfs}]
\btriangle|rba|/<-`<-`<-^)/<1000,600>[\QSup`\QCat`\ \QSymCat;\FP_{\sfc}`\FP_{\sfs}`]
\place(-70,350)[\mbox{\rotatebox{-90}{$\bot$}}]
\place(480,-80)[\bot]
\place(540,390)[\mbox{\rotatebox{150}{$\bot$}}]
\efig
\end{equation}

This question is of crucial importance in the study of \emph{quantale-valued sets} \cite{Hoehle1991,Hoehle1992,Hoehle1995b,Hoehle1998,Hoehle2005,Hoehle2011a}. Let
\[\sQ=(\sQ,\with,k,{}^{\circ})\]
be an \emph{involutive quantale} \cite{Mulvey1992}, considered as the table of truth values. Recall that a \emph{$\sQ$-set} \cite{Hoehle2011a} is exactly a \emph{symmetric} category enriched in a \emph{quantaloid} $\HsQ$ constructed from $\sQ$ (cf. Definition \ref{Q-Set-def} and \cite[Proposition 6.3]{Hoehle2011a}). The category
\[\QSet:=\HsQ\text{-}{\bf SymCat}\]
of $\sQ$-sets is precisely the category of symmetric $\HsQ$-categories. So, the following question becomes a special case of Question \ref{ques}:

\begin{ques} \label{what-is-Q-powerset}
What is the \emph{$\sQ$-powerset} of a $\sQ$-set?
\end{ques}

The main result of this paper, Theorem \ref{s-monadic}, gives an affirmative answer to Question \ref{ques}. Therefore, all the three forgetful functors in the diagram
\[\bfig
\btriangle|llb|/->`->`<-^)/<800,500>[\QSup`\QCat`\ \QSymCat;\FU_{\sfc}`\Us`]
\morphism(800,0)/<-^)/<800,0>[\ \QSymCat`\ \Set/\CQ_0;]
\morphism(0,500)/@{->}@<3pt>/<1600,-500>[\QSup`\ \Set/\CQ_0;\FU]
\efig\]
are strictly monadic, and consequently:
\begin{itemize}
\item The ``power'' of a symmetric $\CQ$-category $X$ is given by its image under $\Us\FP_{\sfs}$ (where $\FP_{\sfs}\dv\Us$), which is precisely the \emph{symmetrization} of the presheaf $\CQ$-category of $X$.
\end{itemize}

In particular, for an involutive quantale $\sQ$, the monad generated by the adjunction $\FP_{\sfs}\dv\Us$ provides an explicit answer to Question \ref{what-is-Q-powerset}; that is,
\begin{itemize}
\item The ``$\sQ$-powerset'' of a $\sQ$-set $X$ is given by its image under $\Us\FP_{\sfs}$ (where $\FP_{\sfs}\dv\Us$), which is precisely the \emph{symmetrization} of the presheaf $\HsQ$-category of $X$.
\end{itemize}
Therefore,  the \emph{$\sQ$-powerset monad} on $\QSet$ is precisely formulated, and we elaborate the details of its components in Section \ref{Ps-Monad-QSet}.

\section{Categories and symmetric categories enriched in a quantaloid}

Complete lattices and join-preserving maps constitute a symmetric monoidal closed category $\Sup$ \cite{Joyal1984}. A \emph{quantaloid} \cite{Rosenthal1996} $\CQ$ is a category enriched in $\Sup$; that is, a category whose hom-sets are complete lattices, such that the composition of $\CQ$-arrows preserves joins on both sides, i.e.,
\[v\circ\Big(\bv_{i\in I}u_i\Big)=\bv_{i\in I}v\circ u_i\quad\text{and}\quad \Big(\bv_{i\in I}v_i\Big)\circ u=\bv_{i\in I}v_i\circ u\]
for all $\CQ$-arrows $u,u_i\colon p\to q$, $v,v_i\colon q\to r$ $(i\in I)$. The corresponding right adjoints induced by the compositions 
\[(-\circ u)\dv(-\lda u)\colon\CQ(p,r)\to\CQ(q,r)\quad\text{and}\quad(v\circ -)\dv(v\rda -)\colon\CQ(p,r)\to\CQ(p,q)\]
satisfy 
\[v\circ u\leq w\iff v\leq w\lda u\iff u\leq v\rda w\] 
for all $\CQ$-arrows $u\colon p\to q$, $v\colon q\to r$, $w\colon p\to r$, where $\lda$ and $\rda$ are called \emph{left} and \emph{right implications} in $\CQ$, respectively.

A \emph{homomorphism} of quantaloids is a functor of the underlying categories that preserves joins of morphisms. A quantaloid $\CQ$ is \emph{involutive} if it is equipped with an \emph{involution}; that is, a homomorphism
\begin{equation} \label{involution}
(-)^{\circ}\colon\CQ^{\op}\to\CQ
\end{equation}
of quantaloids with
\[q^{\circ}=q\quad\text{and}\quad u^{\circ\circ}=u\]
for all $q\in\CQ_0(=\ob\CQ)$ and $\CQ$-arrows $u\colon p\to q$, which necessarily satisfies
\[(1_q)^{\circ}=1_q,\quad (v\circ u)^{\circ}=u^{\circ}\circ v^{\circ}\quad\text{and}\quad\Big(\bv_{i\in I}u_i\Big)^{\circ}=\bv_{i\in I}u_i^{\circ}\]
for all $q\in\ob\CQ$ and $\CQ$-arrows $u,u_i\colon p\to q$, $v\colon q\to r$ $(i\in I)$. 

Throughout this paper, we let $\CQ$ denote a \emph{small} involutive quantaloid; that is, $\CQ$ has a set $\CQ_0$ of objects, and $\CQ$ is equipped with an involution \eqref{involution}.

A \emph{$\CQ_0$-typed set} $X$ is a set $X$ equipped with a \emph{type} map $|\text{-}|:X\to\CQ_0$. With \emph{type-preserving} maps as morphisms, i.e., maps $f\colon X\to Y$ satisfying $|x|=|fx|$ for all $x\in X$, we obtain a category 
\[\Set/\CQ_0.\]

A \emph{$\CQ$-relation} $\phi\colon X\rto Y$ between $\CQ_0$-typed sets consists of a family of $\CQ$-arrows $\phi(x,y)\in\CQ(|x|,|y|)$ $(x\in X,\ y\in Y)$. $\CQ_0$-typed sets and $\CQ$-relations constitute a (\emph{not} necessarily involutive!) quantaloid $\QRel$, in which
\begin{itemize}
\item the local order is inherited from $\CQ$, i.e.,
\[\phi\leq\psi\colon X\rto Y\iff\forall x\in X,\ \forall y\in Y\colon\phi(x,y)\leq\psi(x,y);\]
\item the composition and implications of $\CQ$-relations $\phi\colon X\rto Y$, $\psi\colon Y\rto Z$, $\eta\colon X\rto Z$ are given by
\begin{align*}
&\psi\circ\phi\colon X\rto Z,\quad(\psi\circ\phi)(x,z)=\bv_{y\in Y}\psi(y,z)\circ\phi(x,y),\\
&\eta\lda\phi\colon Y\rto Z,\quad(\eta\lda\phi)(y,z)=\bw_{x\in X}\eta(x,z)\lda\phi(x,y),\\
&\psi\rda\eta\colon X\rto Z,\quad(\psi\rda\eta)(x,y)=\bw_{y\in Y}\psi(y,z)\rda\eta(x,z);
\end{align*}
\item the identity $\CQ$-relation on a $\CQ_0$-typed set $X$ is given by 
\[\id_X\colon X\rto X,\quad \id_X(x,y)=\begin{cases}
1_{|x|} & \text{if}\ x=y,\\
\bot_{|x|,|y|} & \text{else,}
\end{cases}.\]
where $\bot_{|x|,|y|}$ refers to the bottom $\CQ$-arrow in $\CQ(|x|,|y|)$.
\end{itemize}

A \emph{$\CQ$-category} \cite{Rosenthal1996,Stubbe2005,Stubbe2014} consists of a $\CQ_0$-typed set $X$ and a $\CQ$-relation $\al\colon X\rto X$ such that $\id_X\leq\al$ and $\al\circ\al\leq\al$; that is, 
\[1_{|x|}\leq\al(x,x)\quad\text{and}\quad\al(y,z)\circ\al(x,y)\leq\al(x,z)\]
for all $x,y,z\in X$. The underlying (pre)order of a $\CQ$-category $(X,\al)$ is given by
\[x\leq y\iff |x|=|y|\ \text{and}\ 1_{|x|}\leq\al(x,y).\]
We write $x\cong y$ if $x\leq y$ and $y\leq x$. A $\CQ$-category $(X,\al)$ is \emph{separated} (also \emph{skeletal}) if $x=y$ whenever $x\cong y$ in its underlying order.

A \emph{$\CQ$-functor} (resp. \emph{fully faithful $\CQ$-functor}) $f:(X,\al)\to(Y,\be)$ between $\CQ$-categories is a type-preserving map $f\colon X\to Y$ such that
\[\al(x,x')\leq\be(fx,fx')\quad(\text{resp.}\ \al(x,x')=\be(fx,fx'))\]
for all $x,x'\in X$. With the pointwise (pre)order between $\CQ$-functors given by
\[f\leq g\colon (X,\al)\to(Y,\be)\iff\forall x\in X\colon fx\leq gx\iff\forall x\in X\colon 1_{|x|}\leq\be(fx,gx),\]
$\CQ$-categories and $\CQ$-functors constitute a locally ordered category
\[\QCat.\]

A pair of $\CQ$-functors $f\colon(X,\al)\to(Y,\be)$ and $g\colon(Y,\be)\to(X,\al)$ forms an adjunction in $\QCat$, denoted by $f\dv g$, if
\[1_X\leq gf\quad\text{and}\quad fg\leq 1_Y,\]
or equivalently, if
\[\be(fx,y)=\al(x,gy)\]
for all $x\in X$, $y\in Y$. In this case, $f$ is called a \emph{left adjoint} of $g$, and $g$ is a \emph{right adjoint} of $f$.

A $\CQ$-category $(X,\al)$ is \emph{symmetric} \cite{Heymans2011} if
\begin{equation} \label{sym-Q-cat}
\al(x,y)=\al(y,x)^{\circ}
\end{equation}
for all $x,y\in X$. The full subcategory of $\QCat$ consisting of symmetric $\CQ$-categories is denoted by 
\[\QSymCat.\]
From each $\CQ$-category $(X,\al)$ we may construct a symmetric $\CQ$-category $(X,\al_{\sfs})$, with
\begin{equation} \label{symmetrization-def}
\al_{\sfs}(x,y)=\al(x,y)\wedge\al(y,x)^{\circ}
\end{equation}
for all $x,y\in X$. It is clear that $f\colon(X,\al_{\sfs})\to(Y,\be_{\sfs})$ is a $\CQ$-functor whenever so is $f\colon(X,\al)\to(Y,\be)$, giving rise to the \emph{symmetrization} functor
\begin{equation} \label{s-functor}
(-)_{\sfs}\colon\QCat\to\QSymCat.
\end{equation}
In fact, $\QSymCat$ is a coreflective subcategory of $\QCat$, with $(-)_{\sfs}$ being the coreflector \cite{Heymans2011}:

\begin{lem} \label{QSymCat-coref-QCat}
Let $(X,\al)$, $(Y,\be)$ be $\CQ$-categories. If $(X,\al)$ is symmetric, then $f\colon(X,\al)\to(Y,\be)$ is a $\CQ$-functor if, and only if, $f\colon(X,\al)\to(Y,\be_{\sfs})$ is a $\CQ$-functor.
\end{lem}

For $\CQ$-categories $(X,\al)$, $(Y,\be)$, a $\CQ$-relation $\phi\colon X\rto Y$ becomes a \emph{$\CQ$-distributor} $\phi\colon(X,\al)\oto(Y,\be)$ if 
\[\be\circ\phi\circ\al\leq\phi;\] 
that is,
\[\be(y,y')\circ\phi(x,y)\circ\al(x',x)\leq\phi(x',y')\]
for all $x,x'\in X$, $y,y'\in Y$.
$\CQ$-categories and $\CQ$-distributors constitute a quantaloid $\QDist$ which includes $\QRel$ as a full subquantaloid. Compositions and implications of $\CQ$-distributors are computed in the same way as in $\QRel$,  and the identity $\CQ$-distributor on a $\CQ$-category $(X,\al)$ is given by $\al\colon(X,\al)\oto(X,\al)$.

Each $\CQ$-functor $f\colon(X,\al)\to(Y,\be)$ induces an adjunction $f_{\nat}\dv f^{\nat}$ in $\QDist$ (i.e., $\al\leq f^{\nat}\circ f_{\nat}$ and $f_{\nat}\circ f^{\nat}\leq\be$), given by
\[f_{\nat}\colon(X,\al)\oto(Y,\be),\quad f_{\nat}(x,y)=\be(fx,y)\quad\text{and}\quad f^{\nat}\colon(Y,\be)\oto(X,\al),\quad f^{\nat}(y,x)=\be(y,fx),\]
called the \emph{graph} and \emph{cograph} of $f$, respectively. Obviously, the identity $\CQ$-distributor $\al$ is the cograph of the identity $\CQ$-functor $1_X\colon(X,\al)\to(X,\al)$. Hence, if no confusion arises, in what follows we write
\[1_X^{\nat}=\al\]
for the hom of a $\CQ$-category $X=(X,\al)$, and write $X_{\sfs}$ for the symmetrization of $X$. 




For each $q\in\ob\CQ$, let $\{q\}$ denote the (necessarily symmetric) one-object $\CQ$-category whose only object has type $q$ and hom $1_q$. A \emph{presheaf} $\mu$ (of type $q$) on a $\CQ$-category $X$ is a $\CQ$-distributor $\mu\colon X\oto\{q\}$, and presheaves on $X$ constitute a separated $\CQ$-category $\PX$ with
\[1_{\PX}^{\nat}(\mu,\mu')=\mu'\lda\mu\]
for all $\mu,\mu'\in\PX$. Dually, the separated $\CQ$-category $\PdX$ of \emph{copresheaves} on $X$ consists of $\CQ$-distributors $\lam\colon\{q\}\oto X$ with $|\lam|=q$ and
\[1_{\PdX}^{\nat}(\lam,\lam')=\lam'\rda\lam\]
for all $\lam,\lam'\in\PdX$. In particular, for each $q\in\CQ_0$, $\sP\{q\}$ (resp. $\sPd\{q\}$) consists of $\CQ$-arrows of domain (resp. codomain) $q$ as objects.

For every $\CQ$-functor $f\colon X\to Y$, it is straightforward to check that 
\begin{equation} \label{fra-def}
f^{\ra}\colon\PX\to\PY,\quad f^{\ra}\mu=\mu\circ f^{\nat}\quad\text{and}\quad f^{\la}\colon\PY\to\PX,\quad f^{\la}\lam=\lam\circ f_{\nat}
\end{equation}
define an adjunction $f^{\ra}\dv f^{\la}$ in $\QCat$.


\section{Complete categories enriched in a quantaloid}

A $\CQ$-category $X$ is \emph{complete} if the \emph{Yoneda embedding}
\[\sy_X\colon X\to\PX,\quad x\mapsto 1_X^{\nat}(-,x)\]
admits a left adjoint $\sup_X\colon\PX\to X$ in $\QCat$, which is equivalent to the existence of a right adjoint $\inf_X\colon\PdX\to X$ to the \emph{co-Yoneda embedding}
\[\syd_X\colon X\to\PdX,\quad x\mapsto 1_X^{\nat}(x,-).\]

\begin{rem} \label{PXs-sup}
Let $X$ be a complete $\CQ$-category. Elaborating the adjunction $\sup_X\dv\sy_X$ in details, we obtain that
\begin{equation} \label{sup-def}
1_X^{\nat}({\sup}_X\mu,-)=1_X^{\nat}\lda\mu
\end{equation}
for all $\mu\in\PX$. In fact, even if a $\CQ$-relation $\mu\colon X\rto q$ $(q\in\CQ_0)$ is not a presheaf on $X$, its supremum $\sup_X\mu$ still exists, and it is an object of $X$ satisfying \eqref{sup-def}. To see this, just note that $\mu\circ 1_X^{\nat}\in\PX$, and
\[1_X^{\nat}({\sup}_X(\mu\circ 1_X^{\nat}),-)=1_X^{\nat}\lda(\mu\circ 1_X^{\nat})=(1_X^{\nat}\lda 1_X^{\nat})\lda\mu=1_X^{\nat}\lda\mu;\]
that is, $\sup_X\mu={\sup}_X(\mu\circ 1_X^{\nat})$. In particular, the $\CQ$-functor $\sup_X\colon\PX\to X$ can be extended to
\[{\sup}_X\colon\PX_{\sfs}\to X,\]
since a presheaf on $X_{\sfs}$ is always a $\CQ$-relation with domain $X$.
\end{rem}

\begin{exmp} \label{PX-complete}
For each $\CQ$-category $X$, both $\PX$ and $\PdX$ are separated complete $\CQ$-categories. In particular, for each $\Phi\in\sP\PX$ (see \cite[Example 2.9]{Shen2013a}),
\begin{equation} \label{sup-PX}
{\sup}_{\PX}\Phi=\Phi\circ(\sy_X)_{\nat}=\bv_{\mu\in\PX}\Phi(\mu)\circ\mu.
\end{equation}
\end{exmp}

In a $\CQ$-category $X$, the \emph{tensor} of $u\in\sP\{|x|\}$ and $x\in X$, denoted by $u\otimes x$, is an object of $X$ of type $|u\otimes x|=|u|$, such that
\begin{equation} \label{tensor-def}
1_X^{\nat}(u\otimes x,-)=1_X^{\nat}(x,-)\lda u.
\end{equation}
$X$ is \emph{tensored} if $u\otimes x$ exists for all choices of $u$ and $x$. The dual notions are \emph{cotensors} and \emph{cotensored} $\CQ$-categories.


A $\CQ$-category $X$ is \emph{order-complete} if, for any $q\in\CQ_0$,
\[X_q:=\{x\in X\mid |x|=q\}\]
admits all joins (or equivalently, all meets) in its underlying order. In particular, if $X$ is separated and order-complete, then each $X_q$ is a complete lattice.

\begin{prop} \label{QCat-complete-tensored} (See \cite{Stubbe2006}.)
A $\CQ$-category is complete if, and only if, it is tensored, cotensored and order-complete.
\end{prop}

Let $X$ be a complete $\CQ$-category and $x\in X$, $q\in\CQ_0$. For each subset $\{x_i\mid i\in I\}\subseteq X_q$, it follows from \cite[Proposition 3.5.4]{Shen2014} that
\begin{equation} \label{X-x-bw-xi}
1_X^{\nat}\Big(\bv_{i\in I}x_i,x\Big)=\bw_{i\in I}1_X^{\nat}(x_i,x)\quad\text{and}\quad 1_X^{\nat}\Big(x,\bw_{i\in I}x_i\Big)=\bw_{i\in I}1_X^{\nat}(x,x_i),
\end{equation}
where the joins and meets of $x_i$ $(i\in I)$ are computed in the underlying order of $X$. As a consequence, we deduce the following lemma that will be useful later:

\begin{lem} \label{tensor-join}
Let $X$ be a complete $\CQ$-category and $q\in\CQ_0$, $u\in\sP\{q\}$. Then for each subset $\{x_i\mid i\in I\}\subseteq X_q$, 
\[u\otimes\Big(\bv\limits_{i\in I}x_i\Big)=\bv\limits_{i\in I}u\otimes x_i,\]
where the joins are computed in the underlying order of $X$.
\end{lem}

\begin{proof}
Note that 
\begin{align*}
1_X^{\nat}\Big(\bv_{i\in I}u\otimes x_i,-\Big)&=\bw_{i\in I}1_X^{\nat}(u\otimes x_i,-) & \text{(Equations \eqref{X-x-bw-xi})}\\
&=\bw_{i\in I}(1_X^{\nat}(x_i,-)\lda u) & \text{(Equation \eqref{tensor-def})}\\
&=\Big(\bw_{i\in I}(1_X^{\nat}(x_i,-)\Big)\lda u\\
&=1_X^{\nat}\Big(\bv_{i\in I}x_i,-\Big)\lda u. & \text{(Equations \eqref{X-x-bw-xi})}
\end{align*}
The conclusion thus follows from the definition \eqref{tensor-def} of tensors.
\end{proof}

When the $\CQ$-categories under concern are tensored, the $\CQ$-functoriality of a type-preserving map between them can be characterized as follows:

\begin{prop} \label{Qfunctor-tensored} (See \cite{Shen2014}.)
A type-preserving map $f\colon X\to Y$ between tensored $\CQ$-categories is a $\CQ$-functor if, and only if,
\begin{enumerate}[label={\rm(\arabic*)}]
\item \label{Qfunctor-tensored:tensor} $u\otimes_Y fx\leq f(u\otimes_X x)$ for all $x\in X$, $u\in\sP\{|x|\}$, and 
\item \label{Qfunctor-tensored:order} $f$ is an order-preserving map between the underlying ordered sets of $X$, $Y$.
\end{enumerate}
\end{prop}


Furthermore, left adjoint $\CQ$-functors between complete $\CQ$-categories have the following equivalent characterizations:

\begin{prop} (See \cite{Stubbe2005,Stubbe2006}.) \label{left-adjoint-preserves-sup}
For a $\CQ$-functor $f\colon X\to Y$ between complete $\CQ$-categories, the following statements are equivalent:
\begin{enumerate}[label={\rm(\roman*)}]
\item $f$ is a left adjoint in $\QCat$.
\item $f$ is a left adjoint between the underlying ordered sets of $X$, $Y$, and preserves tensors in the sense that $f(u\otimes_X x)=u\otimes_Y fx$ for all $x\in X$, $u\in\sP\{|x|\}$.
\item $f$ is \emph{$\sup$-preserving} in the sense that $f\sup_X=\sup_Y f^{\ra}$. 
\end{enumerate}
\end{prop}

Separated complete $\CQ$-categories and left adjoint $\CQ$-functors (or equivalently, $\sup$-preserving $\CQ$-functors) constitute a subcategory of $\QCat$, and we denote it by
\[\QSup.\]
It is well known (see, e.g., \cite[Proposition 6.11]{Stubbe2005}) that the forgetful functor $\FU_{\sfc}\colon\QSup\to\QCat$ admits a left adjoint
\[\FP_{\sfc}\colon\QCat\to\QSup,\]
which sends each $\CQ$-functor $f\colon X\to Y$ to the left adjoint $\sQ$-functor (see \eqref{fra-def})
\[f^{\ra}\colon\PX\to\PY.\]

Since Lemma \ref{QSymCat-coref-QCat} implies that the inclusion functor $\QSymCat\ \to/^(->/\QCat$ is left adjoint to $(-)_{\sfs}$, it follows soon that the functor
\[\FP_{\sfs}:=(\QSymCat\ \to/^(->/\QCat\to^{\FP_{\sfc}}\QSup)\]
is left adjoint to
\[\Us:=(\QSup\to^{\FU_{\sfc}}\QCat\to^{(-)_{\sfs}}\QSymCat),\]
whose unit and counit are given by
\[\{\sy_X\colon X\to(\PX)_{\sfs}\}_{X\in\QSymCat}\quad\text{and}\quad\{{\sup}_X\colon\PX_{\sfs}\to X\}_{X\in\QSup},\]
respectively, where 
\begin{itemize}
\item $\sy_X\colon X\to(\PX)_{\sfs}$ is the symmetrization (see \eqref{s-functor}) of the Yoneda embedding $\sy_X\colon X\to\PX$, and
\item ${\sup}_X\colon\PX_{\sfs}\to X$ is the extension of ${\sup}_X\colon\PX\to X$ (see Remark \ref{PXs-sup}).
\end{itemize}
The induced monad on $\QSymCat$ is denoted by
\begin{equation} \label{bbPs-def}
\bbPs=(\Ps,\sy,{\sup}_{\sP}),
\end{equation}
where 
\[\Ps:=(\QSymCat\ \to/^(->/\QCat\to^{\FP_{\sfc}}\QSup\to^{\Us}\QSymCat)\]
sends each symmetric $\CQ$-category $X$ to $(\PX)_{\sfs}$.

\begin{rem} \label{Ps-unit-m}
The unit of the monad $\bbPs$ is simply $\{\sy_X\colon X\to(\PX)_{\sfs}\}_{X\in\QSymCat}$. To understand the multiplication
\[{\sup}_{\PX}\colon\Ps\Ps X=(\sP(\PX)_{\sfs})_{\sfs}\to\Ps X=(\PX)_{\sfs},\]
just note that it is the symmetrization (see \eqref{s-functor}) of the $\CQ$-functor
\begin{equation} \label{sup-PXs}
{\sup}_{\PX}\colon\sP(\PX)_{\sfs}\to\PX.
\end{equation}
By Remark \ref{PXs-sup}, \eqref{sup-PXs} is the extension of the $\CQ$-functor 
\[{\sup}_{\PX}\colon\sP\PX\to\PX\]
described by \eqref{sup-PX} in Example \ref{PX-complete}. Indeed, the supremum of each $\Phi\in\sP(\PX)_{\sfs}$ is precisely
\[{\sup}_{\PX}\Phi={\sup}_{\PX}(\Phi\circ 1_{\PX}^{\nat})=\Phi\circ 1_{\PX}^{\nat}\circ(\sy_X)_{\nat}=\Phi\circ (\sy_X)_{\nat}=\bv_{\mu\in\PX}\Phi(\mu)\circ\mu,\]
where $\sy_X$ refers to the original Yoneda embedding $X\to\PX$ as in Example \ref{PX-complete}.
\end{rem}

The purpose of this section is to show that the right adjoint functor $\Us\colon\QSup\to\QSymCat$ is strictly monadic. To this end, we need some preparations.

A \emph{$\CQ$-closure operator} on a $\CQ$-category $X$ is a $\CQ$-functor $c\colon X\to X$ such that
\[1_X\leq c\quad\text{and}\quad cc\cong c.\]
It is easy to see that each pair of adjoint $\CQ$-functors $f\dv g\colon Y\to X$ gives rise to a $\CQ$-closure operator $gf\colon X\to X$. Moreover:

\begin{prop} \label{closure-operator-complete} (See \cite{Shen2013a}.)
Let $c\colon X\to X$ be a $\CQ$-closure operator, and let
\[\Fix(c):=\{x\in X\mid cx\cong x\}\]
be the $\CQ$-subcategory of $X$ consisting of fixed points of $c$.
\begin{enumerate}[label={\rm(\arabic*)}]
\item The inclusion $\CQ$-functor $\Fix(c)\ \to/^(->/X$ is right adjoint to the codomain restriction $\overline{c}\colon X\to\Fix(c)$ of $c$.
\item If $X$ is a complete $\CQ$-category, then so is $\Fix(c)$.
\end{enumerate}
\end{prop}

Recall that in a category $\CC$, an object $B$ is called a \emph{retract} \cite{Riehl2016} of an object $A$ if there are morphisms $f:A\to B$ and $g:B\to A$ such that 
\[fg=1_B.\]
In this case, $f$ is called a \emph{retraction} of $A$ onto $B$, and $g$ is a \emph{section} of $h$. 

\begin{prop} \label{section-lift}
Let $(X,\al)$ be a separated complete $\CQ$-category, and let $Y$ be a retract of $X$ in $\Set/\CQ_0$, with $f\colon X\to Y$ being a retraction. Suppose that
\begin{enumerate}[label={\rm (\alph*)}]
\item \label{section-lift:join} $f\big(\bv\limits_{i\in I}x_i\big)=f\big(\bv\limits_{i\in I}x'_i\big)$ whenever $fx_i=fx'_i$ for all $i\in I$, and
\item \label{section-lift:tensor} $f(u\otimes x)=f(u\otimes x')$ whenever $fx=fx'$ and $u\in\sP\{|x|\}$.
\end{enumerate}
Then there exists a section $h\colon Y\to X$ of $f$ such that
\begin{enumerate}[label={\rm (\arabic*)}]
\item \label{section-lift:complete} $(Y,\be)$ is a separated complete $\CQ$-category with $\be(y,y')=\al(hy,hy')$ for all $y,y'\in Y$,
\item \label{section-lift:dv} $f\dv h\colon(Y,\be)\to(X,\al)$ in $\QCat$, and
\item \label{section-lift:max} if $g\colon Y\to X$ is a section of $f$ in $\Set/\CQ_0$, then $gy\leq hy$ and $\al(gy,gy')\leq\al(hy,hy')$ for all $y,y'\in Y$.
\end{enumerate}
\end{prop}

\begin{proof}
Let $g\colon Y\to X$ be a section of $f$ in $\Set/\CQ_0$. For each $y\in Y$, define
\begin{equation} \label{By-h-def}
B_y:=\{x\in X\mid fx=y\}\quad\text{and}\quad hy:=\bv B_y,
\end{equation}
where the join is computed in the underlying order of the separated complete $\CQ$-category $(X,\al)$, and it is well defined because $gy\in B_y$. Then \ref{section-lift:join} guarantees that 
\begin{equation} \label{fhy=fgy=y}
fhy=fgy=y
\end{equation} 
for all $y\in Y$; that is, $h\colon Y\to X$ is a section of $f$ in $\Set/\CQ_0$.

Let $\be(y,y')=\al(hy,hy')$ for all $y,y'\in Y$. Then $(Y,\be)$ is clearly a $\CQ$-category, which is embedded into $(X,\al)$ via the fully faithful and injective $\CQ$-functor\ $h:(Y,\be)\to(X,\al)$. Next, we show that $f\colon(X,\al)\to(Y,\be)$ is a $\CQ$-functor, which necessarily follows from the $\CQ$-functoriality of $hf\colon(X,\al)\to(X,\al)$. Since $(X,\al)$ is tensored by Proposition \ref{QCat-complete-tensored}, it suffices to check that $hf$ satisfies the two conditions given in Proposition \ref{Qfunctor-tensored}.

First, $hf$ preserves the underlying order of $(X,\al)$. Suppose that $x\leq x'$. Then
\[fx'=f(x\vee x')=f(hfx\vee hfx'),\]
where the second equality follows from $fx=fhfx$, $fx'=fhfx'$ and \ref{section-lift:join}. Hence $hfx\vee hfx'\leq hfx'$ by the definition of $h$, and consequently $hfx\leq hfx'$. 

Second, $u\otimes hfx\leq hf(u\otimes x)$ for all $x\in X$, $u\in\sP\{|x|\}$. Indeed, note that each $z\in B_{fx}$ satisfies $fz=fx$, and consequently $f(u\otimes z)=f(u\otimes x)$ by \ref{section-lift:tensor}; that is, $u\otimes z\in B_{f(u\otimes x)}$. It follows that
\begin{align*}
u\otimes hfx&=u\otimes\bv B_{fx} & \text{(Equations \eqref{By-h-def})}\\
&=\bv\{u\otimes z\mid z\in B_{fx}\} & \text{(Lemma \ref{tensor-join})}\\
&\leq\bv B_{f(u\otimes x)} & (u\otimes z\in B_{f(u\otimes x)}\ \text{if}\ z\in B_{fx})\\
&=hf(u\otimes x). & \text{(Equations \eqref{By-h-def})}
\end{align*}

Therefore, $f\colon(X,\al)\to(Y,\be)$ is a $\CQ$-functor, which already satisfies $fh=1_Y$ by Equation \eqref{fhy=fgy=y}. Since $1_X\leq hf$ is an immediate consequence of the definition of $h$, it follows that 
\[f\dv h\] 
in $\QCat$. In particular, $hf$ is a $\CQ$-closure operator on the separated complete $\CQ$-category $(X,\al)$ and, consequently, the $\CQ$-category $(Y,\be)$ is separated and complete, because it is isomorphic to the $\CQ$-subcategory $\Fix(hf)$ of $(X,\al)$ (see Proposition \ref{closure-operator-complete}).

Finally, for any $y,y'\in Y$, it is clear that $gy\leq hy$. Since $fhy=fgy=y$, it follows from \ref{section-lift:tensor} that
\begin{equation} \label{section-lift:a}
f(\al(gy,gy')\otimes hy)=f(\al(gy,gy')\otimes gy).
\end{equation}
Note that 
\[1_{|y'|}\leq\al(gy,gy')\lda\al(gy,gy')=\al(\al(gy,gy')\otimes gy,gy')\]
implies that $\al(gy,gy')\otimes gy\leq gy'$, which in combination with \eqref{section-lift:a} gives rise to
\[\al(gy,gy')\otimes hy\leq hf(\al(gy,gy')\otimes hy)=hf(\al(gy,gy')\otimes gy)\leq hfgy'=hy'.\]
Hence 
\[1_{|y'|}\leq\al(\al(gy,gy')\otimes hy,hy')=\al(hy,hy')\lda\al(gy,gy');\]
that is, $\al(gy,gy')\leq\al(hy,hy')$.
\end{proof}

Recall that given a functor $\sG\colon\CD\to\CC$:
\begin{itemize}
\item A \emph{G-split coequalizer} is a pair $X\two^f_g Y$ of $\CD$-morphisms such that $\sG X\two^{\sG f}_{\sG g}\sG Y$ extends to a \emph{split coequalizer} diagram
\begin{equation} \label{split-coequalizer}
\bfig
\morphism(0,0)|a|/@{->}@<3pt>/<600,0>[\sG X`\sG Y;\sG f]
\morphism(0,0)|b|/@{->}@<-3pt>/<600,0>[\sG X`\sG Y;\sG g]
\morphism(600,0)|b|/{@{->}@/^1.6em/}/<-600,0>[\sG Y`\sG X;t]
\morphism(600,0)|a|/->/<600,0>[\sG Y`Z;h]
\morphism(1200,0)|b|/{@{->}@/^1em/}/<-600,0>[Z`\sG Y;s]
\efig
\end{equation}
in $\CC$, which means that
\begin{equation} \label{split-def}
h(\sG f)=h(\sG g),\quad hs=1_Z,\quad (\sG g)t=1_{\sG Y}\quad\text{and}\quad (\sG f)t=sh.
\end{equation}
\item $\sG$ \emph{strictly creates coequalizers of $\sG$-split pairs} if, for every $\sG$-split coequalizer \eqref{split-coequalizer}, there exists a unique $\CD$-object $W$ and a unique $\CD$-morphism $k\colon Y\to W$ such that $\sG W=Z$, $\sG k=h$ and 
\[\bfig
\morphism(0,0)|a|/@{->}@<3pt>/<600,0>[X`Y;f]
\morphism(0,0)|b|/@{->}@<-3pt>/<600,0>[X`Y;g]
\morphism(600,0)|a|/->/<600,0>[Y`W;k]
\efig\]
is a coequalizer diagram.
\item A right adjoint functor $\sG\colon\CD\to\CC$ is \emph{strictly monadic} if the canonical comparison functor from $\CD$ to the Eilenberg-Moore category of the induced monad on $\CC$ defines an isomorphism of categories (see, e.g., \cite[Section II.3.2]{Hofmann2014} and \cite[Section 5.3]{Riehl2016}, for details).
\item Beck's monadicity theorem \cite{Beck1967,MacLane1998,Riehl2016} states that a right adjoint functor $\sG\colon\CD\to\CC$ is strictly monadic if, and only if, it strictly creates coequalizers of $\sG$-split pairs (see, e.g., \cite[Theorem 5.5.1 and Exercise 5.5.i]{Riehl2016}). 
\end{itemize}

\begin{thm} \label{s-monadic}
The right adjoint functor $\Us\colon\QSup\to\QSymCat$ is strictly monadic.
\end{thm}

\begin{proof}
It suffices to show that $\Us\colon\QSup\to\QSymCat$ strictly creates coequalizers of $\Us$-split pairs. Let 
\[\bfig
\morphism(0,0)|a|/@{->}@<3pt>/<900,0>[(X,\al)`(Y,\be);f]
\morphism(0,0)|b|/@{->}@<-3pt>/<900,0>[(X,\al)`(Y,\be);g]
\efig\]
be a pair of left adjoint $\CQ$-functors between separated complete $\CQ$-categories such that 
\[\bfig
\morphism(0,0)|a|/@{->}@<3pt>/<900,0>[(X,\al_{\sfs})`(Y,\be_{\sfs});f]
\morphism(0,0)|b|/@{->}@<-3pt>/<900,0>[(X,\al_{\sfs})`(Y,\be_{\sfs});g]
\morphism(900,0)|b|/{@{->}@/^1.6em/}/<-900,0>[(Y,\be_{\sfs})`(X,\al_{\sfs});t]
\morphism(900,0)|a|/->/<900,0>[(Y,\be_{\sfs})`(Z,\ga);h]
\morphism(1800,0)|b|/{@{->}@/^1em/}/<-900,0>[(Z,\ga)`(Y,\be_{\sfs});s]
\efig\]
is a split coequalizer diagram in $\QSymCat$, which by Equations \eqref{split-def} means that
\begin{equation} \label{split-QCat}
hf=hg,\quad hs=1_Z,\quad gt=1_Y\quad\text{and}\quad ft=sh.
\end{equation}

{\bf Step 1.} $h\colon Y\to Z$ satisfies the conditions of Proposition \ref{section-lift}, which induces a section $s'\colon Z\to Y$ such that 
\begin{enumerate}[label=(\arabic*)]
\item \label{sprime:complete} $\xi(z,z')=\be(s'z,s'z')$ defines a separated complete $\CQ$-category $(Z,\xi)$,
\item \label{sprime:dv} $h\dv s'\colon(Z,\xi)\to(Y,\be)$ in $\QCat$, and
\item \label{sprime:max} $sz\leq s'z$ and $\be(sz,sz')\leq\be(s'z,s'z')$ for all $z,z'\in Z$.
\end{enumerate}
Moreover, $\ga=\xi_{\sfs}$.

First, $h\big(\bv\limits_{i\in I}y_i\big)=h\big(\bv\limits_{i\in I}y'_i\big)$ whenever $hy_i=hy'_i$ for all $i\in I$. In this case, it follows from \eqref{split-QCat} that 
\[fty_i=shy_i=shy'_i=fty'_i.\]
The combination of Proposition \ref{left-adjoint-preserves-sup} and \eqref{split-QCat} then implies that
\begin{align*}
h\big(\bv\limits_{i\in I}y_i\big)&=h\big(\bv\limits_{i\in I}gty_i\big)=hg\big(\bv\limits_{i\in I}ty_i\big)=hf\big(\bv\limits_{i\in I}ty_i\big)=h\big(\bv\limits_{i\in I}fty_i\big)\\
&=h\big(\bv\limits_{i\in I}fty'_i\big)=hf\big(\bv\limits_{i\in I}ty'_i\big)=hg\big(\bv\limits_{i\in I}ty'_i\big)=h\big(\bv\limits_{i\in I}gty'_i\big)=h\big(\bv\limits_{i\in I}y'_i\big).
\end{align*}

Second, $h(u\otimes y)=h(u\otimes y')$ whenever $hy=hy'$ and $u\in\sP\{|y|\}$. In this case, by applying \eqref{split-QCat} and Proposition \ref{left-adjoint-preserves-sup} again, we deduce that 
\[fty=shy=shy'=fty',\] 
and consequently
\begin{align*}
h(u\otimes y)&=h(u\otimes gty)=hg(u\otimes ty)=hf(u\otimes ty)=h(u\otimes fty)\\
&=h(u\otimes fty')=hf(u\otimes ty')=hg(u\otimes ty')=h(u\otimes gty')=h(u\otimes y').
\end{align*}
as desired. 

Finally, $\ga=\xi_{\sfs}$. Let $z,z'\in Z$. On one hand, since $\ga$ is symmetric, from the functoriality of $s$ and \ref{sprime:max} we obtain that
\[\ga(z,z')\leq\be(sz,sz')\wedge\be(sz',sz)^{\circ}\leq\be(s'z,s'z')\wedge\be(s'z',s'z)^{\circ}=\xi_{\sfs}(z,z').\]
On the other hand, 
\[\xi_{\sfs}(z,z')=\be_{\sfs}(s'z,s'z')\leq\ga(hs'z,hs'z')=\ga(z,z').\]

{\bf Step 2.} $\bfig
\morphism(0,0)|a|/@{->}@<3pt>/<600,0>[(X,\al)`(Y,\be);f]
\morphism(0,0)|b|/@{->}@<-3pt>/<600,0>[(X,\al)`(Y,\be);g]
\morphism(600,0)|a|/->/<600,0>[(Y,\be)`(Z,\xi);h]
\efig$ is a coequalizer diagram in $\QSup$. Let $h'\colon(Y,\be)\to(Z',\xi')$ be a left adjoint $\CQ$-functor between separated complete $\CQ$-categories satisfying $h'f=h'g$. We claim that $h's'\colon(Z,\xi)\to(Z',\xi')$ is the unique $\CQ$-functor that makes the right triangle of the diagram
\[\bfig
\morphism(0,0)|a|/@{->}@<3pt>/<900,0>[(X,\al)`(Y,\be);f]
\morphism(0,0)|b|/@{->}@<-3pt>/<900,0>[(X,\al)`(Y,\be);g]
\morphism(900,0)|a|/->/<900,0>[(Y,\be)`(Z,\xi);h]
\morphism(900,0)|l|/->/<900,-400>[(Y,\be)`(Z',\xi');h']
\morphism(1800,0)|r|/-->/<0,-400>[(Z,\xi)`(Z',\xi');h's']
\efig\]
commutative. On one hand, note that for any $y,y'\in Y$, if $hy=hy'$, then
\[h'y=h'gty=h'fty=h'shy=h'shy'=h'fty'=h'gty'=h'y',\]
which in conjunction with Proposition \ref{left-adjoint-preserves-sup} implies that
\[h's'hy=h'\Big(\bv\{y'\in Y\mid hy'=hy\}\Big)=\bv\{h'y'\mid y'\in Y,\ hy'=hy\}=h'y\]
for all $y\in Y$; that is, the right triangle of the above diagram is commutative. On the other hand, if $h''\colon(Z,\xi)\to(Z',\xi')$ satisfies $h''h=h'$, then
\[h''=h''hs'=h's'.\]
It remains to show that $h's'\colon(Z,\ga)\to(Z',\ga')$ is a left adjoint in $\QCat$. To this end, note that $h'$ has a right adjoint $t'\colon(Z',\xi')\to(Y,\be)$ in $\QCat$. Since
\[h's'ht'=h't'\leq 1_{Z'}\quad\text{and}\quad ht'h's'\geq hs'=1_Z,\]
we conclude that $h's'\dv ht'$, as desired.

{\bf Step 3.} For the uniqueness of the lifting of $(Z,\ga)$ to a separated complete $\CQ$-category, suppose that $(Z,\eta)$ is another separated complete $\CQ$-category such that 
\[\bfig
\morphism(0,0)|a|/@{->}@<3pt>/<600,0>[(X,\al)`(Y,\be);f]
\morphism(0,0)|b|/@{->}@<-3pt>/<600,0>[(X,\al)`(Y,\be);g]
\morphism(600,0)|a|/->/<600,0>[(Y,\be)`(Z,\eta);h] 
\efig\]
is a coequalizer diagram in $\QSup$. Then there exists a unique left adjoint $\CQ$-functor $k\colon(Z,\eta)\to(Z,\xi)$ that makes the right triangle of the diagram
\[\bfig
\morphism(0,0)|a|/@{->}@<3pt>/<900,0>[(X,\al)`(Y,\be);f]
\morphism(0,0)|b|/@{->}@<-3pt>/<900,0>[(X,\al)`(Y,\be);g]
\morphism(900,0)|a|/->/<900,0>[(Y,\be)`(Z,\eta);h]
\morphism(900,0)|l|/->/<900,-400>[(Y,\be)`(Z,\xi);h]
\morphism(1800,0)|r|/-->/<0,-400>[(Z,\eta)`(Z,\xi);k]
\efig\]
commutative; that is, $kh=h$. Thus, by \eqref{split-QCat} it is easy to see that
\[kz=khsz=hsz=z\]
for all $z\in Z$, which forces $k=1_Z$. So, the identity map $1_Z\colon(Z,\eta)\to(Z,\xi)$ is a left adjoint $\CQ$-functor, whose right adjoint must be given by $1_Z\colon(Z,\xi)\to(Z,\eta)$. Hence, the $\CQ$-functoriality of $1_Z$ on both sides forces $\xi(z,z')=\eta(z,z')$ for all $z,z'\in Z$, which completes the proof.
\end{proof}

\begin{cor}
The Eilenberg-Moore category $\QSymCat^{\bbPs}$ is isomorphic to $\QSup$. Hence, $\QSup$ is strictly monadic over $\QSymCat$.
\end{cor}

\section{The powerset monad on quantale-valued sets} \label{Ps-Monad-QSet}

A \emph{(unital) quantale} \cite{Mulvey1986,Rosenthal1990} is exactly a one-object quantaloid. Throughout this section, we let
\[\sQ=(\sQ,\with,k,{}^{\circ})\]
denote an involutive quantale. Explicitly:
\begin{itemize}
\item $\sQ$ is a complete lattice (with a top element $\top$ and a bottom element $\bot$).
\item $(\sQ,\with,k)$ is a monoid, such that the multiplication $\with$ preserves joins on both sides.
\item The left and right implications and induced by the multiplication are denoted by $\ldd$ and $\rdd$, respectively, which satisfy
\[p\with q\leq r\iff p\leq r\ldd q\iff q\leq p\rdd r\]
for all $p,q,r\in\sQ$.
\item $\sQ$ is equipped with an involution, i.e., a map $(-)^{\circ}\colon\sQ\to\sQ$ such that
\[k^{\circ}=k,\quad q^{\circ\circ}=q,\quad (p\with q)^{\circ}=q^{\circ}\with p^{\circ}\quad\text{and}\quad\Big(\bv_{i\in I}q_i\Big)^{\circ}=\bv_{i\in I}q_i^{\circ}\]
for all $p,q,q_i\in\sQ$.
\end{itemize}

From $\sQ$ we may construct a quantaloid $\HsQ$ \cite{Hoehle2011a}, given by the following data:
\begin{itemize}
\item Objects of $\HsQ$ are \emph{hermitian} (also \emph{self-adjoint}) elements of $\sQ$; that is, $q\in\sQ$ satisfying $q^{\circ}=q$.
\item Given hermitian elements $p,q\in\sQ$, $\HsQ(p,q)$ consists of elements $d\in\sQ$ satisfying
\begin{equation} \label{HsQ-mor-def}
d\leq p\wedge q\quad\text{and}\quad(d\ldd p)\with p=d=q\with(q\rdd d).
\end{equation}
\item The composition of $d\in\HsQ(p,q)$ and $e\in\HsQ(q,r)$ is given by 
\begin{equation} \label{HsQ-comp-def}
e\circ d:=(e\ldd q)\with d=e\with(q\rdd d).
\end{equation}
\item The identity morphism on $q\in\sQ$ is $q$ itself.
\item Each hom-set $\HsQ(p,q)$ is equipped with the order inherited from $\sQ$.
\end{itemize}



$\HsQ$ is obviously an involutive quantaloid with the involution lifted from $\sQ$. From the definition we see that a $\HsQ$-category consists of a set $X$, a map $|\text{-}|:X\to\sQ$ and a map $\al\colon X\times X\to\sQ$ such that
\begin{enumerate}[label=(\arabic*)]
\item \label{HsQ-cat:str} $\al(x,y)\leq|x|\wedge|y|$,
\item \label{HsQ-cat:div} $(\al(x,y)\ldd |x|)\with |x|=\al(x,y)=|y|\with(|y|\rdd\al(x,y))$,
\item \label{HsQ-cat:ref} $|x|\leq\al(x,x)$,
\item \label{HsQ-cat:tran} $(\al(y,z)\ldd |y|)\with\al(x,y)=\al(y,z)\with(|y|\rdd\al(x,y))\leq\al(x,z)$
\end{enumerate}
for all $x,y,z\in X$, where \ref{HsQ-cat:str} and \ref{HsQ-cat:div} follows from $\al(x,y)\in\HsQ(|x|,|y|)$. Note that the combination of \ref{HsQ-cat:str} and \ref{HsQ-cat:ref} forces
\[\al(x,x)=|x|\]
for all $x\in X$, and thus a $\HsQ$-category is exactly given by a map $\al:X\times X\to\sQ$ such that 
\begin{enumerate}[label=(S\arabic*)]
\item \label{sim-def:str} $\al(x,y)\leq\al(x,x)\wedge\al(y,y)$,
\item \label{sim-def:div} $(\al(x,y)\ldd\al(x,x))\with\al(x,x)=\al(x,y)=\al(y,y)\with(\al(y,y)\rdd\al(x,y))$,
\item \label{sim-def:tran} $(\al(y,z)\ldd\al(y,y))\with\al(x,y)=\al(y,z)\with(\al(y,y)\rdd\al(x,y))\leq\al(x,z)$
\end{enumerate}
for all $x,y,z\in X$, and it is symmetric if
\begin{enumerate}[label=(S\arabic*),start=4]
\item \label{sim-def:sym} $\al(x,y)=\al(y,x)^{\circ}$
\end{enumerate}
for all $x,y\in X$.

\begin{defn} \label{Q-Set-def} (See \cite{Hoehle2011a}.)
A \emph{$\sQ$-set} is a symmetric $\HsQ$-category; that is, a set $X$ equipped with a map $\al:X\times X\to\sQ$ satisfying \ref{sim-def:str}--\ref{sim-def:sym}.
\end{defn}

\begin{rem}
The notion of $\HsQ$ here is slightly different from \cite{Lai2020}. The quantaloid $\HsQ$ in \cite{Lai2020} has all elements of $\sQ$ as its objects, while in this paper we restrict the objects of $\HsQ$ to hermitian elements of $\sQ$. Nevertheless, as \cite[Remark 4.1]{Lai2020} reveals, it makes no difference when we only deal with symmetric $\HsQ$-categories.
\end{rem}

\begin{rem} \label{QSet-discussion}
A $\sQ$-set may be viewed as a set $X$ equipped with a \emph{$\sQ$-valued equality} (or \emph{$\sQ$-valued similarity}) $\al$ \cite{Hoehle2011a,Lai2020}. The value $\al(x,y)$ is interpreted as the extent of $x$ being equal to $y$, and $\al(x,x)$ represents the extent of existence of $x$ (since every entity is supposed to be equal to itself). Therefore:
\begin{itemize}
\item \ref{sim-def:str} says that $x$ is equal to $y$ only if both $x$ and $y$ exist.
\item The first equality of \ref{sim-def:div} says that $x$ is equal to $y$ if, and only if, $x$ exists and its existence forces $x$ being equal to $y$.
\item The first inequality of \ref{sim-def:tran} says that if $x$ is equal to $y$, and the existence of $y$ forces $y$ being equal to $z$, then $x$ is equal to $z$. 
\item \ref{sim-def:sym} says that if $x$ is equal to $y$, then $y$ is equal to $x$.
\end{itemize}
\end{rem}

\begin{exmp} \label{Q-set-exmp}
Some important examples of $\sQ$-sets are listed below:
\begin{enumerate}[label=(\arabic*)]
\item \label{Q-set-exmp:2} If $\sQ={\bf 2}$, the two-element Boolean algebra, then a ${\bf 2}$-set $(X,\al)$ is just an equivalence relation on a subset of $X$. Explicitly,
\[\{(x,y)\in X\times X\mid\al(x,y)=1\}\]
is an equivalence relation on the subset $\{x\in X\mid\al(x,x)=1\}$ of $X$, whose elements are supposed to ``exist''. In particular, $(X,\al)$ reduces to a (crisp) set if
\begin{itemize}
\item $(X,\al)$ is separated, i.e., $\al(x,y)=1$ if and only if $x=y$;
\item $(X,\al)$ is \emph{global}, i.e., $\al(x,x)=1$ for all $x\in X$. 
\end{itemize}
\item \label{Q-set-exmp:PCX} If $\sQ$ is a \emph{frame}, then $\sQ$-sets are precisely \emph{$\Om$-sets} in the sense of Fourman-Scott \cite{Fourman1979}. In particular, given a topological space $X$, let
\[\PCX:=\{f\mid f\ \text{is a real-valued continuous map on an open subset}\ D(f)\subseteq X\}.\]
For any $f,g\in\PCX$, let
\[\al(f,g):=\Int\{x\in D(f)\cap D(g)\mid f(x)=g(x)\},\]
i.e., the interior of the subset of $X$ consisting of elements on which $f$ and $g$ coincide. Then $(\PCX,\al)$ is an $\OX$-set, where $\OX$ is the frame of open subsets of $X$.
\item \label{Q-set-exmp:metric} Let $\sQ$ be the Lawvere quantale $[0,\infty]=([0,\infty],+,0)$ \cite{Lawvere1973}. Then $[0,\infty]$-sets are symmetric \emph{partial metric spaces}; that is, sets $X$ equipped with a map 
\[\al\colon X\times X\to[0,\infty]\]
such that
\[\al(x,x)\vee\al(y,y)\leq\al(x,y),\quad\al(x,z)\leq\al(y,z)-\al(y,y)+\al(x,y)\quad\text{and}\quad\al(x,y)=\al(y,x)\]
for all $x,y,z\in X$. In particular, let
\[\CI:=\{[a,b]\mid 0\leq a<b\leq \infty\}\]
be the set of closed intervals contained in $[0,\infty]$. Then
\[\al([a,b],[c,d])=b\vee d-a\wedge c\]
defines a symmetric partial metric space $(\CI,\al)$.
\end{enumerate}
\end{exmp}

We denote by
\[\QSet:=\HsQ\text{-}{\bf SymCat}\]
the category of $\sQ$-sets, whose morphisms are maps $f\colon(X,\al)\to(Y,\be)$ between $\sQ$-sets satisfying
\begin{equation} \label{QSet-f}
\al(x,x)=\be(fx,fx)\quad\text{and}\quad\al(x,x')\leq\be(fx,fx')
\end{equation}
for all $x,x'\in X$. Following the interpretations of Remark \ref{QSet-discussion}, \eqref{QSet-f} says that $x$ exists if and only if $fx$ exists, and if $x$ is equal to $x'$, then $fx$ is equal to $fx'$.

Now, let us elaborate how the monad
\begin{equation} \label{bbPs-QSet} 
\bbPs=(\Ps,\sy,{\sup}_{\sP})
\end{equation}
given by \eqref{bbPs-def} on the category $\HsQ\text{-}{\bf SymCat}$ describes the \emph{$\sQ$-powerset monad} on $\QSet$.

First of all, for each $\sQ$-set $(X,\al)$,
\[\Ps(X,\al)\]
is the \emph{$\sQ$-powerset} of $(X,\al)$, whose elements are $\HsQ$-distributors $\mu\colon(X,\al)\oto\{q\}$ $(q\in\sQ)$; that is, maps $\mu\colon X\to\sQ$ such that
\begin{enumerate}[label=(P\arabic*)]
\item \label{PsX:str} $\mu(x)\leq\al(x,x)\wedge q$,
\item \label{PsX:div} $(\mu(x)\ldd\al(x,x))\with\al(x,x)=\mu(x)=q\with(q\rdd\mu(x))$,
\item \label{PsX:dist} $(\mu(y)\ldd\al(y,y))\with\al(x,y)=\mu(y)\with(\al(y,y)\rdd\al(x,y))\leq\mu(x)$
\end{enumerate}
for all $x,y\in X$. 

\begin{defn} \label{potential-subset}
A \emph{potential $\sQ$-subset} of a $\sQ$-set $(X,\al)$ is a pair $(\mu,q)$, where $\mu\colon X\to\sQ$ and $q\in\sQ$ satisfies \ref{PsX:str}--\ref{PsX:dist}.
\end{defn}

So, the $\sQ$-powerset of a $\sQ$-set consists of its potential $\sQ$-subsets, which can be understood as follows:

\begin{rem} \label{potential-subset-disucssion}
In a potential $\sQ$-subset $(\mu,q)$ of a $\sQ$-set $(X,\al)$:
\begin{itemize}
\item the value $\mu(x)$ represents the degree of $x$ being in $(\mu,q)$, and 
\item $q$ represents the degree of $(\mu,q)$ being a $\sQ$-subset of $(X,\al)$. 
\end{itemize}
Therefore:
\begin{itemize}
\item \ref{PsX:str} says that $x$ is in $(\mu,q)$ only if $x$ exists and $(\mu,q)$ is a subset of $(X,\al)$.
\item The first equality of \ref{PsX:div} says that $x$ is in $(\mu,q)$ if, and only if, $x$ exists and its existence forces $x$ being in $(\mu,q)$. The second equality of \ref{PsX:div} says that $x$ is in $(\mu,q)$ if, and only if, $(\mu,q)$ is a subset of $(X,\al)$ and this fact forces $x$ being in $(\mu,q)$.
\item The first inequality of \ref{PsX:dist} says that if $x$ is equal to $y$, and the existence of $y$ forces $y$ being in $(\mu,q)$, then $x$ is in $(\mu,q)$.
\end{itemize}
\end{rem}
 
\begin{exmp} \label{potential-Q-subset-exmp}
For the examples listed in \ref{Q-set-exmp}:
\begin{enumerate}[label=(\arabic*)]
\item A potential ${\bf 2}$-subset of a ${\bf 2}$-set $(X,\al)$ is either $(\varnothing,0)$ or $(U,1)$, where $U$ is a subset of $A:=\{x\in X\mid\al(x,x)=1\}$ that is a union of some equivalence classes of the corresponding equivalence relation on $A$; in other words, if  $y\in U$ and $x$ is equivalent to $y$, then $x\in U$. In particular, if $(X,\al)$ is separated, then $U$ can be any subset of $A$.
\item Let $(\PCX,\al)$ be the $\OX$-set considered in Example \ref{Q-set-exmp}\ref{Q-set-exmp:PCX}. A potential $\OX$-subset of $(\PCX,\al)$ is a pair $(\mu,V)$, where $\mu$ is a map $\mu\colon\PCX\to\OX$ such that
\[\al(f,g)\cap\Int(\mu(g)\cup(X\setminus D(g)))=\mu(g)\cap\Int(\al(f,g)\cup(X\setminus D(g)))\subseteq\mu(f)\subseteq D(f)\cap V\]
for all $f,g\in\PCX$, where $X\setminus D(g)$ refers to the complement of the set $D(g)$ in $X$.
\item Let $(X,\al)$ be a symmetric partial metric space (see Example \ref{Q-set-exmp}\ref{Q-set-exmp:metric}). A potential $[0,\infty]$-subset of $(X,\al)$ is pair $(\mu,q)$, where $\mu$ is a map $\mu\colon X\to[0,\infty]$ such that
\[\al(x,x)\vee q\leq\mu(x)\leq\mu(y)+\al(x,y)-\al(y,y)\]
for all $x,y\in X$. In particular, a potential $[0,\infty]$-subset of the symmetric partial metric space $(\CI,\al)$ is a pair $(\mu,q)$, where $\mu$ is a map $\mu\colon\CI\to[0,\infty]$ such that
\[(b-a)\vee q\leq\mu([a,b])\leq\mu([c,d])+b\vee d-a\wedge c-d+c\]
for all $[a,b],[c,d]\in\CI$.
\end{enumerate}
\end{exmp}

Once the notion of potential $\sQ$-subset is made clear, it is straightforward to interpret the components of the $\sQ$-powerset monad \eqref{bbPs-QSet} as the $\sQ$-valued version of \eqref{fra}, \eqref{y} and \eqref{bigcup}:
\begin{itemize}
\item The functor $\Ps$ sends a map $f\colon(X,\al)\to(Y,\be)$ in $\QSet$ to the map
\[f^{\ra}\colon\Ps(X,\al)\to\Ps(Y,\be)\]
between the corresponding $\sQ$-powersets. Explicitly, for each potential $\sQ$-subset $(\mu,q)$ of $(X,\al)$, 
\[f^{\ra}(\mu,q):=(\lam,q)\]
is a potential $\sQ$-subset of $(Y,\be)$, with
\begin{equation} \label{lam-y-def}
\lam(y)=\bv_{x\in X}(\mu(x)\ldd\al(x,x))\with\be(y,fx)
\end{equation}
for all $y\in Y$. Obviously, \eqref{lam-y-def} says that $y$ is in $(\lam,q)$ if, and only if, there exists $x$ such that $x$ is in $(\mu,q)$ and $y$ is equal to $fx$.
\item The unit of \eqref{bbPs-QSet} is given by
\[\sy_{(X,\al)}\colon(X,\al)\to\Ps(X,\al),\quad x\mapsto(\al(-,x),\al(x,x)),\]
where $(\al(-,x),\al(x,x))$ of $(X,\al)$ is the potential $\sQ$-subset of $(X,\al)$ such that
\begin{itemize}
\item the degree of $(\al(-,x),\al(x,x))$ being a $\sQ$-subset of $(X,\al)$ is the same as the extent of existence of $x$, and
\item $y$ is in $(\al(-,x),\al(x,x))$ if, and only if, $y$ is equal to $x$.
\end{itemize}
\item By Remark \ref{Ps-unit-m}, the multiplication of \eqref{bbPs-QSet} is given by
\[{\sup}_{(X,\al)}\colon\Ps\Ps(X,\al)\to\Ps(X,\al),\quad(\Phi,p)\mapsto\Big(\bv_{(\mu,q)\in\Ps(X,\al)}(\Phi(\mu,q)\ldd q)\with\mu,p\Big),\quad\]
which means that $x$ is in ${\sup}_{(X,\al)}(\Phi,p)$ if, and only if, there exists a potential $\sQ$-subset $(\mu,q)$ of $(X,\al)$ such that $(\mu,q)$ is in $(\Phi,p)$ and $x$ is in $(\mu,q)$.
\end{itemize}

\section*{Acknowledgement}

The first named author acknowledges the support of National Natural Science Foundation of China (No. 12071319). The authors would like to thank Professor Hongliang Lai and Professor Dexue Zhang for helpful discussions.






\end{document}